\documentclass[11pt]{article}

\RequirePackage[amsthm,amsmath]{imsart}

\usepackage[mathscr]{eucal}
\usepackage{graphicx}
\usepackage{latexsym}
\usepackage{amssymb}
\usepackage{psfrag}
\usepackage{epsfig}
\DeclareMathAlphabet{\mathpzc}{OT1}{pzc}{m}{it}

\include{psbox}

\usepackage{amsfonts}
\usepackage{graphicx}
\usepackage{amsfonts}
\usepackage{color}
\usepackage[usenames,dvipsnames]{xcolor}
\usepackage[textsize=small]{todonotes}

\usepackage[numbers]{natbib}

\begin{document}

\newtheorem{theorem}{\bf Theorem}[section]
\newtheorem{example}{\bf Example}[section]
\newtheorem{definition}{\bf Definition}[section]
\newtheorem{corollary}{\bf Corollary}[section]
\newtheorem{remark}{\bf Remark}[section]
\newtheorem{lemma}{\bf Lemma}[section]
\newtheorem{assumption}{Assumption}
\newtheorem{condition}{\bf Condition}[section]
\newtheorem{proposition}{\bf Proposition}[section]
\newtheorem{definitions}{\bf Definition}[section]
\newtheorem{problem}{\bf Problem}
\numberwithin{equation}{section}
\newcommand{\skp}{\vspace{\baselineskip}}
\newcommand{\noi}{\noindent}
\newcommand{\osc}{\mbox{osc}}

\newcommand{\img}{\imath}
\newcommand{\iy}{\infty}
\newcommand{\eps}{\varepsilon}
\newcommand{\del}{\delta}
\newcommand{\Rk}{\mathbb{R}^k}
\newcommand{\R}{\mathbb{R}}
\newcommand{\spa}{\vspace{.2in}}
\newcommand{\V}{\mathcal{V}}
\newcommand{\E}{\mathbb{E}}
\newcommand{\I}{\mathbb{I}}
\newcommand{\p}{\mathbb{P}}
\newcommand{\PP}{\mathbb{P}}
\newcommand{\sgn}{\mbox{sgn}}

\newcommand{\QQ}{\mathbb{Q}}

\newcommand{\lan}{\langle}
\newcommand{\ran}{\rangle}
\newcommand{\lf}{\lfloor}
\newcommand{\rf}{\rfloor}
\def\wh{\widehat}
\newcommand{\defn}{\stackrel{def}{=}}
\newcommand{\txb}{\tau^{\epsilon,x}_{B^c}}
\newcommand{\tyb}{\tau^{\epsilon,y}_{B^c}}
\newcommand{\tilxb}{\tilde{\tau}^\eps_1}
\newcommand{\pxeps}{\mathbb{P}_x^{\eps}}
\newcommand{\non}{\nonumber}
\newcommand{\dist}{\mbox{dist}}

\newcommand{\Om}{\mathnormal{\Omega}}
\newcommand{\om}{\omega}
\newcommand{\vph}{\varphi}
\newcommand{\Del}{\mathnormal{\Delta}}
\newcommand{\Gam}{\mathnormal{\Gamma}}
\newcommand{\Sig}{\mathnormal{\Sigma}}

\newcommand{\tilyb}{\tilde{\tau}^\eps_2}
\newcommand{\beq}{\begin{eqnarray*}}
\newcommand{\eeq}{\end{eqnarray*}}
\newcommand{\beqn}{\begin{eqnarray}}
\newcommand{\eeqn}{\end{eqnarray}}
\newcommand{\ink}{\rule{.5\baselineskip}{.55\baselineskip}}

\newcommand{\bt}{\begin{theorem}}
\newcommand{\et}{\end{theorem}}
\newcommand{\deps}{\Del_{\eps}}

\newcommand{\be}{\begin{equation}}
\newcommand{\ee}{\end{equation}}
\newcommand{\ac}{\mbox{AC}}
\newcommand{\BB}{\mathbb{B}}
\newcommand{\VV}{\mathbb{V}}
\newcommand{\DD}{\mathbb{D}}
\newcommand{\KK}{\mathbb{K}}
\newcommand{\HH}{\mathbb{H}}
\newcommand{\TT}{\mathbb{T}}
\newcommand{\CC}{\mathbb{C}}
\newcommand{\ZZ}{\mathbb{Z}}
\newcommand{\SSS}{\mathbb{S}}
\newcommand{\EE}{\mathbb{E}}

\newcommand{\clg}{\mathcal{G}}
\newcommand{\clb}{\mathcal{B}}
\newcommand{\cls}{\mathcal{S}}
\newcommand{\clc}{\mathcal{C}}
\newcommand{\cld}{\mathcal{D}}
\newcommand{\cle}{\mathcal{E}}
\newcommand{\clv}{\mathcal{V}}
\newcommand{\clu}{\mathcal{U}}
\newcommand{\clr}{\mathcal{R}}
\newcommand{\clt}{\mathcal{T}}
\newcommand{\cll}{\mathcal{L}}

\newcommand{\cli}{\mathcal{I}}
\newcommand{\clp}{\mathcal{P}}
\newcommand{\cla}{\mathcal{A}}
\newcommand{\clf}{\mathcal{F}}
\newcommand{\clh}{\mathcal{H}}
\newcommand{\N}{\mathbb{N}}
\newcommand{\Q}{\mathbb{Q}}
\newcommand{\bfx}{{\boldsymbol{x}}}
\newcommand{\bfb}{{\boldsymbol{b}}}
\newcommand{\bfw}{{\boldsymbol{w}}}
\newcommand{\bfz}{{\boldsymbol{z}}}
\newcommand{\bfu}{{\boldsymbol{u}}}
\newcommand{\bfell}{{\boldsymbol{\ell}}}

\newcommand{\curvz}{{\bf \mathpzc{z}}}
\newcommand{\curvx}{{\bf \mathpzc{x}}}
\newcommand{\curvi}{{\bf \mathpzc{i}}}
\newcommand{\curvs}{{\bf \mathpzc{s}}}

\newcommand{\BM}{\mbox{BM}}

\newcommand{\tac}{\mbox{\scriptsize{AC}}}
\newcommand{\beginsec}{
\setcounter{lemma}{0} \setcounter{theorem}{0}
\setcounter{corollary}{0} \setcounter{definition}{0}
\setcounter{example}{0} \setcounter{proposition}{0}
\setcounter{condition}{0} \setcounter{assumption}{0}
\setcounter{remark}{0} }

\numberwithin{equation}{section} \numberwithin{lemma}{section}

\begin{frontmatter}
\title{On the multi-dimensional skew Brownian motion}

 \runtitle{Skew Brownian motion}

\begin{aug}
\author{Rami Atar\thanks{Research
supported in part by the Israel Science Foundation (Grant 1349/08)
and the Technion fund for promotion of research}
and Amarjit Budhiraja\thanks{Research  supported in part by the National Science Foundation (DMS-1004418, DMS-1016441, DMS-1305120) and the Army Research
 Office (W911NF-10-1-0158)}\\ \ \\
}
\end{aug}

\today

\skp

\begin{abstract}
We provide a new, concise proof of weak existence and uniqueness of solutions to
the stochastic differential equation for the multidimensional skew Brownian motion.
We also present an application to Brownian particles with skew-elastic collisions.

\noi {\bf AMS 2000 subject classifications:} Primary 60J60, 60H10; Secondary 93E15, 60J55.

\noi {\bf Keywords:} skew Brownian motion, SDE with local time, singular diffusion processes.
\end{abstract}

\end{frontmatter}

\section{Introduction}\label{introsec}\beginsec

Let $\Sig$ denote the hyperplane $\{x\in\R^d:x_1=0\}$ in $\R^d$, $d\ge1$, and let a vector
field $b:\Sig\to\R^d$ be given on it.
Consider the stochastic differential equation (SDE), for a process
$X$ taking values in $\R^d$, of the form
\begin{equation}\label{00}
X(t)=x+W(t)+\int_0^tb(X(s))dL(s),\qquad t\ge0,
\end{equation}
where $x\in\R^d$, $W$ is a standard $d$-dimensional Brownian motion and $L$ is the local
time of $X$ at $\Sig$, and
$b$ is bounded and Lipschitz, and satisfies $b_1(x)\in[-1,1]$ for all $x\in\Sig$.
This paper provides a new proof of
weak existence and uniqueness of solutions to \eqref{00}.
A more general equation that allows for a bounded, measurable drift coefficient is
also treated. These results are known by the work of
Takanobu \cite{tak-exi}, \cite{tak-uni}.
Our purpose here is to provide a much shorter proof.
In the one-dimensional case this process is known as {\it skew Brownian motion},
first described by It\^o and McKean \cite{ito-mck}, and further studied
by Walsh \cite{wal} and treated by means of a SDE by Harrison and Shepp \cite{har-she}.
Moreover, the equation can be viewed as an extension of
reflected Brownian motion. Indeed, if $b_1(x)=1$ for all $x\in\Sig$
then the equation describes a reflected Brownian motion in the half space with a reflection
vector field given by $b$,
a case first treated by Anderson and Orey \cite{andore}.
See Lejay \cite{lejay} for a comprehensive survey on skew Brownian motion and SDE involving
local time.

The first to consider multidimensional skew Brownian motion is
Portenko in \cite{por1}, \cite{por2} and \cite{por3}.
The point of view of these papers is
to consider a diffusion for which the drift term has a singularity localized at a
(not necessarily flat) surface.
Existence of solutions to the SDE are provided in \cite{por2}.

Zaitseva \cite{zai} studies an equation similar to \eqref{00} with $b$ for which
$b_1(x)$ does not depend on $x$ and
with an additional term $\sigma(X_s)d\tilde B(L(s))$, and obtains
strong existence and uniqueness.
In the special case $\sigma=0$, the additional term alluded to above disappears, and
one obtains precisely the form \eqref{00} with $b_1$ constant. In this case
the projection of $X$ onto $e_1$ gives a one-dimensional skew Brownian
motion and thus the questions of existence and uniqueness of solutions
are considerably simpler than in the situation studied in the present paper.
Additional works on multidimensional skew Brownian motion
include \={O}shima \cite{osh}, Trutnau \cite{tru} and Anulova \cite{anu}.
The most
relevant work is by Takanobu \cite{tak-exi} and \cite{tak-uni},
where an equation of the form \eqref{00} with general drift and diffusion coefficients
is treated.  The first paper treats existence of weak solutions and the second proves uniqueness.
Existence is shown by constructing a suitable approximating diffusion with a non-singular drift
while uniqueness relies on rather involved arguments
from Brownian excursion theory \cite{wata}.
While the results of the present article
are a special case of those established by Takanobu, the proof (of both existence and
uniqueness) provided here is much shorter and more elementary.

Existence of solutions is established by constructing a `skew random walk' and showing that in diffusion scaling it converges in distribution to a weak solution of the equation.
The proof of weak uniqueness provided here is inspired by the technique
used in \cite{wein} in a one-dimensional setting with time-varying skewness.
Specifically, \cite{wein} considers the equation
\begin{equation}
	\label{eq:wein}
	X(t)=x+W(t)+\int_0^t\alpha(s) dL(s),\qquad t\ge0,
\end{equation}
where $d=1$ and $\alpha: [0, \infty) \to [-1,1]$ is a measurable function.  The key idea
is to compute the conditional
characteristic function $E(e^{\img \lambda X(T)}\mid \clf_t)$, for $0 \le t \le T < \infty$, $\lambda \in \R$ where
$\clf_t = \sigma\{X(s), W(s), L(s): 0 \le s \le t\}$.
In the setting of \cite{wein} this calculation reduces to characterizing the conditional law of $L(T)-L(t)$ given $\clf_t$, which turns out to be simply the law of the local time
(at time instant $T-t$) of a one-dimensional reflected Brownian motion.
This is relevant to the multidimensional setting of \eqref{00} if one regards \eqref{eq:wein}
as an equation for the first component $X_1$, and replaces the time-dependent skewness
$\alpha$ by the state-dependent coefficient $b_1(X)$.
Implementing this is, however, somewhat more involved than in \cite{wein}.
In particular, the calculation of the conditional characteristic function requires a
characterization of the {\em joint}
law of $(X_2(T)-X_2(t), L(T)-L(t))$ given $\clf_t$, where we write $X(t) = (X_1(t), X_2(t)) \in \R \times \R^{d-1}$.  Lemmas \ref{lem:skorprob} and \ref{lem:charcndlaw} provide
this characterization, by which the argument can be completed.
We note that \cite{wein} in fact proves
pathwise uniqueness.  The proof of this result relies on the observation that
if $X$ and $\tilde X$ are two solutions of \eqref{eq:wein}
with associated local times $L$ and $\tilde L$ and common Brownian motion $W$
then $X \vee \tilde X$ is also a solution, an idea that goes back to \cite{lg-thesis}
and \cite{perkins}. However, this approach fails in our setting.

The following notation will be used.  For a Polish space $\cls$, $\BM(\cls)$ will denote the space of real bounded measurable functions on $\cls$, $\clb(\cls)$ the Borel $\sigma$-field on $\cls$ and $\clp(\cls)$ will denote the space of probability measures on $\cls$ equipped with the topology of weak convergence.
Also $\CC_{\cls} = \clc([0,\infty): \cls)$ (resp. $\DD_{\cls} = \cld([0,\infty): \cls)$) will denote  the space of continuous functions (resp. RCLL functions) from $[0,\infty)$ to $\cls$
equipped with the usual local uniform topology (resp. Skorohod topology).  A sequence of $\DD_{\cls}$-valued random variables will be said to be $C$-tight if the corresponding sequence of probability laws on $\DD_{\cls}$ is relatively compact and any weak limit point of the sequence is supported on $\CC_{\cls}$.  For $g \in \BM(\cls)$, $\|g\|_{\infty} = \sup_{x \in \cls} |g(x)|$.  Elements in $\R^d$ will be regarded as row vectors.

The rest of this paper is organized as follows. In the next section we state the main result
and present some open problems related to it.
The proof is provided in Section \ref{sec:mainproof}.
Finally, Section \ref{sec4} presents an application to a recent model
of \cite{FIK} for Brownian particles undergoing skew-elastic collisions.

\section{Main result}\label{prelim}
\beginsec

Let $(\Om, \clf, \clp , \{\clf_t\}_{t\ge 0})$
be a filtered probability space and let $\{W(t)\}_{0\le t < \infty}$ be a standard $d$-dimensional $\{\clf_t\}$-Brownian motion.
Given a continuous $\{\clf_t\}$-semimartingale $\{Z(t)\}_{t\ge 0}$, define an
$\{\clf_t\}$-adapted continuous process $\{L^{Z}(t)\}_{t\ge 0}$ as
\begin{equation}\label{01}
L^{Z}(t) = |Z(t)| - |Z(0)|  - \int_0^t \sgn(Z(s)) dZ(s),\quad t \ge 0,
\end{equation}
where for $x \in \R$, $\sgn(x) =  1_{(0, \infty)}(x)  - 1_{(-\infty , 0)}(x)$.
$L^Z$ is a nondecreasing process and is referred to as the symmetric local time of $Z$ at $0$ (see \cite[Chapter VI]{rev-yor}).

We consider a stochastic differential equation for a multidimensional process
that involves its local time at the surface $\Sig=\{x\in\R^d:x_1=0\}$.
The equation takes the form
\begin{equation}
	X(t) = x + W(t) + \int_0^t a(X(s)) ds + \int_0^t b(X(s)) dL(s), \quad t \ge 0,
	\label{eq:maineq}
\end{equation}
where $X = (X_1, \ldots, X_d)$, $x \in \R^d$, $a: \R^d \to \R^d$ is bounded and measurable, $b:\Sig\to \R^d$ is a bounded Lipschitz function satisfying $b_1(x) \in [-1,1]$
for all $x \in \Sig$ where $b_1$ denotes the first coordinate of $b$,
and $L$ is the symmetric local time of $X_1$ at $0$.
Our goal is to show that there exists a unique weak solution of \eqref{eq:maineq}.
We recall below the definition of a weak solution.  Let $\CC_{\R^d} = C([0,\infty): \R^d)$, $\CC_{\R} = C([0,\infty): \R)$ and $\bar \CC = \CC_{\R^d} \times \CC_{\R^d} \times \CC_{\R}$.
Denote by $\{\clg_t\}_{t\ge 0}$ the canonical filtration on $\bar \CC$ and introduce canonical processes $\bfx, \bfw, \bfell$ defined by
\begin{equation}\label{eq:coordproc}
\bfw(t)(\bar \omega) = \omega_1(t),\; \bfx(t)(\bar \omega) = \omega_2(t),\; \bfell(t)(\bar \omega) = \omega_3, \; \bar \omega = (\omega_1, \omega_2, \omega_3),\quad t \ge 0.
\end{equation}
\begin{definition}
	\label{def:weaksolnmain}
	\begin{description}
		\item{(a)} A probability measure $P \in \clp(\bar \CC)$ is said to be  a weak solution of \eqref{eq:maineq} if
		\begin{description}
			\item{(i)} $P$-a.s.,
			$$\bfx(t) = x + \bfw(t) + \int_0^t a(\bfx(s)) ds + \int_0^t b(\bfx(s)) d \bfell(s), \quad t \ge 0,$$
			\item{(ii)} Under $P$, $\bfw$ is a $\{\clg_t\}$-Brownian motion,
			\item{(iii)} Under $P$, $\bfell$ is the symmetric local time of $\bfx$.
		\end{description}
		\item{(b)} We say that \eqref{eq:maineq} has a unique weak solution if whenever $P_1, P_2 \in \clp(\bar \CC)$ are two weak solutions of \eqref{eq:maineq}, one has
		$P_1 =
		P_2$.
	\end{description}
\end{definition}
The following is our main result.
\begin{theorem}
	\label{thm:mainthm}
	There exists a unique weak solution of \eqref{eq:maineq}.
\end{theorem}

We end this section with a discussion of open problems addressing natural
extensions of the above result.

\noi{\bf Problem 1a.}
{\it  Determine whether pathwise uniqueness holds for \eqref{eq:maineq}.}

Note that an affirmative solution of the above would automatically lead to
strong existence.

Next, it is natural to consider a setting where $\Sig$ is replaced by
a smooth surface of co-dimension 1, $\hat\Sig$.
In fact, this setting goes back to the aforementioned papers by Portenko.
Let such a surface by given along with a vector field $b:\hat\Sig\to\R^d$,
and assume
  \[
  b(x)\cdot n(x)\in[-1,1] \text{ for every } x\in\hat\Sig,
  \]
where $n$ denotes a unit normal to $\hat\Sig$. Then equation \eqref{eq:maineq},
with $L$ denoting the local time of $X$ at $\hat\Sig$, is indeed more general than the
one treated in this paper
(one defines local time here by means of limits of occupation time in the vicinity of the
surface as e.g.\ in \cite{doku}).
Note that this equation includes as a special case
Brownian motion reflected on a smooth surface, i.e.\ when $b\cdot n=1$
on $\hat\Sig$, for which existence and uniqueness is well-understood
\cite{andore,lio-szn}.
One may address this equation
by flattening $\hat\Sig$ to $\Sig$, that is, by applying a twice continuously differentiable,
invertible map $\vph:\R^d\to\R^d$ that maps (locally) $\hat\Sig$ into $\Sig$.
Setting $Y=\vph(X)$ one obtains an equation for $Y$ of the form
\begin{equation}\label{02}
Y(t)=y+\int_0^t\sigma(Y(s))dW(s)+\int_0^ta(Y(s))ds+\int_0^tc(Y(s))dL(s),
\end{equation}
where $L$ denotes the local time of $Y$ at $\Sig$ and
the coefficients $\sigma$, $a$ and $c$ depend on $b$ and on derivatives,
up to second order, of $\vph$.
This translates the equation into one that involves the local time on a flat surface,
and SDE with coefficients, and so weak existence and uniqueness are available
by \cite{tak-exi} and \cite{tak-uni}. However, strong form of uniqueness is missing,
which leads us to the following extended version of Problem 1a.

\noi{\bf Problem 1b.}
{\it
Determine whether pathwise uniqueness holds for \eqref{02} for a smooth surface $\hat\Sig$.
}

Finally, consider equation \eqref{eq:maineq} in dimension 2 with discontinuous coefficient
$b=e_1$ on $\{x_1=0, x_2\le0\}$
and $b=-e_1$ on $\{x_1=0,x_2>0\}$ (say, with $a=0$). The authors are not aware of
uniqueness results
for this nonstandard, yet simple equation of a reflected Brownian motion.
It would be interesting to understand this equation in the following broader context.

\noi{\bf Problem 2.}
{\it
Determine whether any form of uniqueness holds for \eqref{eq:maineq}
beyond Lipschitz continuity of the coefficient $b$.
}

\section{Proof}\label{sec:mainproof}
\beginsec

The proof has three parts. We begin by showing that the problem can be reduced to
the case where $a=0$. Namely, we show, in Subsection \ref{s31}
that Theorem \ref{thm:mainthm} follows
from
\begin{proposition}
	\label{prop:zerodrift}
	Suppose that $a=0$. Then \eqref{eq:maineq} has a unique weak solution.
\end{proposition}
\noi
Then we turn to proving Proposition \ref{prop:zerodrift}. Setting $a=0$,
we prove uniqueness and then existence in Subsection \ref{sec:uniqaeq0} and
\ref{sec:exiseq0}, respectively.

\subsection{Reduction to the case $a=0$}\label{s31}

Theorem \ref{thm:mainthm} follows from Proposition \ref{prop:zerodrift} by a straightforward application of Girsanov's theorem, as we now show.

\noindent {\bf Proof of Theorem \ref{thm:mainthm}.}
Let $P_0$ be the unique solution of \eqref{eq:maineq} with $a = 0$.  Define $\tilde P \in \clp(\bar \CC)$ as
$$
\tilde P(A) = \int_A \cle_T dP_0, \; A \in \clg_T, \; T > 0,$$
where
\[
\cle_T = \exp \left\{ \int_0^T a(\bfx(s)) \cdot d\bfw(s)
- \frac{1}{2} \int_0^T |a(\bfx(s))|^2 ds \right\}.
\]
Let $\tilde \bfw(t) = \bfw(t) - \int_0^t a(\bfx(s)) ds$, $t\ge 0$ and let $P \in \clp(\bar \CC)$ be defined
as $P = \tilde P \circ (\tilde \bfw, \bfx, \bfell)^{-1}$.  Then by Girsanov's theorem $P$ is a weak solution of \eqref{eq:maineq}.

To see uniqueness of a weak solution, let $P_1$ be another weak solution of \eqref{eq:maineq}.  Define
$$
\tilde P_1(A) = \int_A \tilde \cle_T dP_1, \; A \in \clg_T, \; T > 0,$$
where
\[
\tilde \cle_T = \exp \left\{ -\int_0^T a(\bfx(s)) \cdot d\bfw(s)
- \frac{1}{2} \int_0^T |a(\bfx(s))|^2 ds \right\}.
\]
Let $\tilde \bfw_1(t) = \bfw(t) + \int_0^t a(\bfx(s)) ds$.  Then
$\tilde P_0 = \tilde P_1 \circ (\tilde \bfw_1, \bfx, \bfell)^{-1}$ is a weak solution of \eqref{eq:maineq} with $a = 0$ and so by Proposition \ref{prop:zerodrift}
$\tilde P_0 = P_0$.
Finally for all $T > 0$ and $f \in \BM(\bar \CC)$
\begin{align*}
	&\hspace{-2em}
    E_1 \left[f\left(\bfw(\cdot \wedge T), \bfx(\cdot \wedge T), \bfell(\cdot \wedge T\right))\right]\\
    &=
	\tilde E_1 \left [ f\left(\bfw(\cdot \wedge T), \bfx(\cdot \wedge T), \bfell(\cdot \wedge T)\right) \tilde \cle_T^{-1}\right]\\
	&= \tilde E_0 \left [ f\left(\bfw(\cdot \wedge T) - \int_0^{\cdot \wedge T} a(\bfx(s)) ds, \bfx(\cdot \wedge T), \bfell(\cdot \wedge T)\right)  \cle_T\right]\\
	&= E_0 \left [ f\left(\bfw(\cdot \wedge T) - \int_0^{\cdot \wedge T} a(\bfx(s)) ds, \bfx(\cdot \wedge T), \bfell(\cdot \wedge T)\right)  \cle_T\right]\\
	&= \tilde E \left[f\left(\bfw(\cdot \wedge T) - \int_0^{\cdot \wedge T} a(\bfx(s)) ds, \bfx(\cdot \wedge T), \bfell(\cdot \wedge T)\right)\right]\\
	&= E \left[f\left(\bfw(\cdot \wedge T), \bfx(\cdot \wedge T), \bfell(\cdot \wedge T)\right)\right].\\
\end{align*}
Thus we have $P_1 = P$.
\qed

For the rest of the section we assume that $a = 0$.

\subsection{Weak uniqueness}
\label{sec:uniqaeq0}

For $x=(x_1,\ldots,x_d)\in\R^d$ write
\[
y=(y_1,y_2)\in\R^d \quad \text{where}\quad
 y_1=|x_1|,\quad y_2=(x_2,\ldots,x_d).
\]
Define $\beta: \R^{d-1} \to \R^d$ as $\beta(z) = b(0, z_1, \ldots, z_{d-1})$, $z = (z_1, \ldots, z_{d-1}) \in \R^{d-1}$.
Write $\beta(z) = (\beta_1(z),\beta_2(z)) \in \R \times \R^{d-1}$, $z \in \R^d$.
Consider the equations
\begin{align}
	S(t) &= y_1 + B_1(t) + V(t), \label{eq:onlyfirst}\\
    Y(t) &= y_2 + B(t) + \int_0^t \beta_2(Y(s)) dV(s). \label{eq:allbutfirst}
\end{align}
Here $B_1$ and $B$ are independent 1- and $(d-1)$-dimensional Brownian motions, respectively,
$S$ is a continuous nonnegative process and $V$ is a continuous nondecreasing
process with $V(0) =0$, which increases only when $S=0$.

The following lemma shows that \eqref{eq:onlyfirst}--\eqref{eq:allbutfirst} has a unique weak solution.  Recall the coordinate processes $\bfw, \bfx, \bfell$ on $\bar \CC$ introduced
in \eqref{eq:coordproc} and write $\bfw = (\bfw_1, \bfw_2)$ and $\bfx = (\bfx_1, \bfx_2)$,
where $\bfw_1$ and $\bfx_1$ take values in $\R$ while
$\bfw_2$ and $\bfx_2$ take values in $\R^{d-1}$.
\begin{lemma}
	\label{lem:skorprob}
		\begin{description}
		\item{(a)} There exists a   $Q \in \clp(\bar \CC)$ such that
		\begin{description}
			\item{(i)}  $Q$-a.s., $\bfx_1(t)\ge0$ for all $t$,
and $t \mapsto \bfell(t)$ is nondecreasing and $\bfell(0) = 0$,
			\item{(ii)} $Q$-a.s.,
			\begin{align*}
				\bfx_1(t) &= y_1 + \bfw_1(t) + \bfell(t),\\
				\bfx_2(t) &= y_2 + \bfw_2(t)+ \int_0^t \beta_2(\bfx_2(s)) d \bfell(s), \quad t \ge 0.
			\end{align*}
			\item{(iii)} Under $Q$, $\bfw_1$ and $\bfw_2$ are independent 1- and
$(d-1)$-dimensional $\{\clg_t\}$-Brownian motions,
			\item{(iv)} $Q$-a.s., $\int 1_{(0,\infty)}(\bfx_1(t)) d\bfell(t) = 0$.
		\end{description}
		\item{(b)} Let $Q_1$ and $Q_2$ be two probability measures as in part (a).  Then $Q_1 = Q_2$.
	\end{description}
\end{lemma}
\noi
{\bf Proof.}
The result is an immediate consequence of the following two observations:
\begin{itemize}
	\item Given a 1-dimensional Brownian motion $B_1$, there is a unique pair of processes $(S, V)$ that satisfies  \eqref{eq:onlyfirst} and is such that $S(t) \ge 0$, $V(0) = 0$, $V$
	is nondecreasing and $\int 1_{(0,\infty)}(S(t)) dV(t) = 0$.  This unique pair is given explicitly as
	$$S(t) = y_1 + B(t) - \inf_{0 \le s \le t} \left \{ y_1 + B(s)\right \} \wedge 0, \; V(t) = (S(t) - y_1 - B_1(t)), \quad t \ge 0.$$
	\item Given a $(d-1)$-dimensional Brownian motion $B$ and $V$ defined as above, we have, using the Lipschitz property of $\beta$ and a standard argument based on Gronwall's lemma, that there
	exists a unique continuous process $Y$ that solves \eqref{eq:allbutfirst}.
\end{itemize}
\qed

In order to indicate the dependence of $Q$ in Lemma \ref{lem:skorprob} on the initial condition $(y_1, y_2) \in \R_+\times \R^{d-1}$ we will denote this probability measure
as $Q_{y_1, y_2}$.  It is easy to see that $(y_1, y_2) \mapsto Q_{y_1,y_2}$ is a continuous function.

Let $P$ be a weak solution of \eqref{eq:maineq} (recall that we are taking $a= 0$)
and for $t \ge 0$ consider the following stochastic processes on $(\bar \CC, \clb(\bar \CC),P)$
\begin{align*}
	\hat \bfb^t(\cdot) &= \int_t^{t+\cdot} \sgn(\bfx_1(s)) d \bfw_1(s), \qquad \hat \bfw_2^t(\cdot) = \bfw_2(t+\cdot) - \bfw_2(t),\\
	\hat \bfx_1^t(\cdot) &= \bfx_1(t+ \cdot),\qquad \hat \bfx_2^t(\cdot) = \bfx_2(t+ \cdot),
\qquad \hat \bfell^t(\cdot) = \bfell(t+ \cdot) - \bfell(t).
\end{align*}
Also let $\bfz = (\bfw, \bfx, \bfell)$, $\bfz^t(\cdot) = \bfz(\cdot \wedge t)$,
$\hat \bfz^t = (\hat \bfb^t, \hat \bfw_2^t, |\hat \bfx_1^t|, \hat \bfx_2^t, \hat \bfell^t)$ and $\bfu^t = (\bfz^t, \hat \bfz^t)$.
Note that, for each $t$,
$\bfu^t$ is a random variable with values in $\bar \CC \times \bar \CC$.  Let $Q^t = P\circ (\bfu^t)^{-1}$ and disintegrate $Q^t$ as
$$
Q^t(d\bar\omega, d\hat\omega) = \clr^t(d\bar \omega) \clp^t(\bar \omega, d\hat \omega), \quad (\bar \omega, \hat \omega) \in \bar \CC \times \bar \CC.$$
Note that $\clr^t$ is the probability law of $\bfz^t$ under $P$ and $\clp^t(\bfz^t, \cdot)$ is the conditional law under $P$ of $\hat \bfz^t$ given $\bfz^t$.

With the above notation we have the following lemma.
\begin{lemma}
	\label{lem:charcndlaw}
	Fix $t \ge 0$. Then for $\clr^t$-a.e.\ $\bar \omega$,
	$\clp^t(\bar \omega, d \hat \omega) = Q_{|\bar \bfx_1(t)|, \bar \bfx_2(t)}(d \hat \omega)$,
	where we write $\bar \om = (\bar \om_1, \bar \om_2, \bar \om_3) \in \CC_{\R^d} \times \CC_{\R^d} \times \CC_{\R}$,
	$(\bar \bfx_1(t), \bar \bfx_2(t)) = \bar \om_2(t) \in \R \times \R^{d-1}$.
\end{lemma}
\noi
{\bf Proof.}
Using Tanaka's formula \eqref{01} and the fact that
$\hat\bfell^t$ is the symmetric local time of $\hat\bfx_1^t$,
(as follows from Definition \ref{def:weaksolnmain}(iii)),
we have
\[
|\hat \bfx_1^t(\cdot)| = |\bfx_1(t)| + \int_0^\cdot\sgn(\hat\bfx_1^t(s))d\hat\bfx_1^t(s)
+ \hat \bfell^t(\cdot).
\]
Using Lemma \ref{lem:skorprob}(v) and the property $\sgn(0)=0$ one has
$\int \sgn(\bfx_1(s))d\bfell(s)=0$. Hence
by the expression for $\bfx_1$ from Definition \ref{def:weaksolnmain}(i), $P$-a.s.,
\begin{align*}
	|\hat \bfx_1^t(\cdot)| &= |\bfx_1(t)| + \hat \bfb^t(\cdot) + \hat \bfell^t(\cdot),\\
	\hat \bfx_2^t(\cdot) &= \bfx_2(t) + \hat \bfw_2^t(\cdot) + \int_0^{\cdot} \beta_2(\hat \bfx_2^t(s)) d \hat \bfell^t(s).
\end{align*}
In view of Lemma \ref{lem:skorprob} it suffices to show that for each fixed $t \ge 0$, under $P$, conditionally on $\bfz^t$,
$\hat \bfw_2^t$ and $\hat \bfb^t$ are independent $(d-1)$-dimensional and 1-dimensional Brownian motions and are adapted to the filtration generated by $\hat \bfz^t$, and
have future increments independent of that filtration. Namely,
$$
P \left [ (\hat \bfb^t, \hat \bfw_2^t) \in A \mid \sigma(\bfz^t)\right] = P(\bfw \in A), \mbox{ for all } A \in \clb(\CC_{\R^d}),$$
for all $f \in \BM(\CC_{\R} \times \CC_{\R^{d-1}})$ and $t_1 \ge 0$
$$
E\left[ f(\hat \bfb^t(\cdot \wedge t_1),\hat \bfw^t_2(\cdot \wedge t_1))  \mid \sigma(\bfz^t) \vee \sigma(\hat \bfz^t(s): 0 \le s \le t_1)\right] =
f(\hat \bfb^t(\cdot \wedge t_1),\hat \bfw^t_2(\cdot \wedge t_1)),$$
and for all $0 \le t_1 \le t_2 < \infty$, $g \in \BM(\bar \CC)$, $h\in\BM(\R^d)$
\begin{align*}
	&E\left[ g(\hat \bfz^t(\cdot \wedge t_1)) h( \hat \bfb^t(t_2) - \hat \bfb^t(t_1), \hat \bfw^t(t_2) - \hat \bfw^t(t_1)) \mid \sigma(\bfz^t)\right] \\
	&\quad =
E\left[ g(\hat \bfz^t(\cdot \wedge t_1)) \mid \sigma(\bfz^t)\right]
E\left[h( \hat \bfb^t(t_2) - \hat \bfb^t(t_1), \hat \bfw^t(t_2) - \hat \bfw^t(t_1)) \mid \sigma(\bfz^t)\right].
\end{align*}
The above properties are immediate from the definitions.
\qed

We can now prove the following.
\begin{proposition}
	\label{prop:atmostonesoln}
    With $a=0$, there is at most one weak solution of \eqref{eq:maineq}.
\end{proposition}
\noi{\bf Proof.}
It suffices to show that for every $m\ge1$
and $0 \le t_1 \le t_2\le \cdots \le t_m < \infty$ the probability law
$P\circ (\bfx(t_1), \ldots, \bfx(t_m))^{-1}$ on $(\R^d)^m$ is the same for any weak solution $P$ of \eqref{eq:maineq}.
For simplicity we only consider the case $m=2$ and write $t_1=t$, $t_2 = T$.
Fix $\lambda \in \R^d$.  Then $\bfz(\cdot) = \exp(\img \lambda \cdot \bfx(\cdot))$ satisfies for all $0 \le t \le s$
$$
\bfz(s) = \bfz(t) + \img \int_t^s \bfz(u) \lambda \cdot d \bfw(u) - \frac{|\lambda|^2}{2} \int_t^s \bfz(u) du + \img \int_t^s \bfz(u) \lambda \cdot \beta(\bfx_2(u)) d\bfell(u).$$
Thus writing for $s \ge t$, $E(\bfz(s) \mid \clg_t) = G(s)$, we have
\begin{equation}\label{eq:Geqn}
G(s) = G(t) - \frac{|\lambda|^2}{2} \int_t^s G(u) du
+ E \left[ \img \int_0^{s-t} \exp(\img \lambda_2 \cdot \hat \bfx_2^t(u)) \lambda \cdot \beta(\hat \bfx_2^t(u)) d \hat \bfell^t(u) \Big| \clg_t\right]\end{equation}
where $\lambda = (\lambda_1, \lambda_2) \in \R \times \R^{d-1}$.
For $(y_1, y_2) \in \R_+ \times \R^{d-1}$ denote the expectation under $Q_{y_1, y_2}$ as $E^Q_{y_1, y_2}$ and let
$$
E^Q_{y_1, y_2} \left [
\img \int_0^{s-t} \exp(\img \lambda_2 \cdot \bfx_2(u)) \lambda \cdot \beta(\bfx_2(u)) d \bfell(u) \right] = H^{t,s}(y_1, y_2).$$
Then from Lemma \ref{lem:charcndlaw} and \eqref{eq:Geqn} we have for all $s \ge t$
$$G(s) = G(t) - \frac{|\lambda|^2}{2} \int_t^s G(u) du  + H^{t,s}(|\bfx_1(t)|, \bfx_2(t)).$$
Thus we have shown that for all $0 \le t \le T < \infty$ and $\lambda \in \R^d$ there is a measurable function $H_{\lambda}^{t,T}$ from $\R_+ \times \R^{d-1}$
to the complex plane such that for every weak solution $P$ of \eqref{eq:maineq}
$$E\left[ \exp(\img \lambda \cdot \bfx(T))\mid \clg_t\right] = H_{\lambda}^{t,T}(|\bfx_1(t)|, \bfx_2(t)),\; P\mbox{-a.s.}$$
By a standard approximation argument it now follows that for all $0 \le t \le T < \infty$ and for every $\varphi \in \BM(\R^d)$ there is a function
$H_{\varphi}^{t,T} \in \BM(\R_+ \times \R^{d-1})$ such that for every weak solution $P$ of \eqref{eq:maineq}
\begin{equation}\label{eq:eqhphi}
	E\left[ \varphi(\bfx(T))\mid \clg_t\right] = H_{\varphi}^{t,T}(|\bfx_1(t)|, \bfx_2(t)),\;
P\mbox{-a.s.}
\end{equation}
Finally let $\varphi_1, \varphi_2 \in \BM(\R^d)$ and let $P$ be a weak solution of \eqref{eq:maineq}.  Then
\begin{align*}
	E\left[\varphi_1(\bfx(t)) \varphi_2(\bfx(T))\right] &= E\left[\varphi_1(\bfx(t)) E \left[ \varphi_2(\bfx(T)) \mid \clg_t \right] \right] \\
	&= E\left[\varphi_1(\bfx(t))H_{\varphi_2}^{t,T}(|\bfx_1(t)|, \bfx_2(t))  \right] \\
\end{align*}
Note that
$$\psi(z) = \varphi_1(z) H_{\varphi_2}^{t,T}(|z_1|, z_2), \; z = (z_1, z_2) \in \R \times \R^{d-1}$$
is in $\BM(\R^d)$.  Thus from \eqref{eq:eqhphi}
$$
E\left[\varphi_1(\bfx(t)) \varphi_2(\bfx(T))\right] = H_{\psi}^{0,t}(|x_1|, x_2).
$$
As a result, for every weak solution $P$ of \eqref{eq:maineq}, the probability law $P\circ (\bfx(t), \bfx(T))^{-1}$ does not depend on $P$.  The result follows. \qed

\subsection{ Weak existence}
\label{sec:exiseq0}

\begin{proposition}
	\label{prop3}
	With $a=0$, there exists a weak solution to \eqref{eq:maineq}.
\end{proposition}

\noi{\bf Proof.}
Recall the function $\beta$ introduced in Section \ref{sec:uniqaeq0}.  Fix $x \in \R^d$. In order to prove
the existence of a weak solution of \eqref{eq:maineq} it suffices
to construct continuous stochastic processes $X,W,L$ taking values in $\R^d$, $\R^d$
and $\R$ respectively on some probability space such that
\begin{itemize}
\item $W$ is a standard $d$-dimensional Brownian motion adapted to the filtration
 \[
 \{\clf_t = \sigma\{X(\cdot \wedge t),W(\cdot \wedge t), L(\cdot \wedge t)\}\}_{t \ge 0}
 \]
    and having increments $W_{t_1}-W_t$ independent of ${\cal F}_t$, $t_1>t$,
\item For all $t \ge 0$
$$X(t) = x + W(t) + \int_0^t \beta(Y(s)) dL(s),$$
where $X(t) = (U(t), Y(t)) \in \R \times \R^{d-1}$,
\item $L$ is the symmetric local time of $U$ at $0$.
\end{itemize}

Write $\beta(\xi) = (\beta_1(\xi), \ldots, \beta_d(\xi))$, $\xi \in \R^{d-1}$.
For $\xi \in \R^{d-1}$ and $i = 2, \ldots, d$, let
$\bar \beta_i(\xi)=2\lfloor (\beta_i(\xi) +1)/2\rfloor$
and $\hat \beta_i(\xi)=\beta_i(\xi)-\bar \beta_i(\xi)$.
Note that $\bar \beta_i$ takes values in $2\ZZ$ and $\hat \beta_i$ in $[-1,1]$.
Let
$$
p_{i,1}(\xi) = \frac{1+\hat \beta_i(\xi)}{2},\quad p_{i,-1}(\xi)=\frac{1-\hat \beta_i(\xi)}{2}, \quad i = 2, \ldots, d, \quad \xi \in \R^{d-1}.$$
Also let
$$
p_{1,1}(\xi) = \frac{1+ \beta_1(\xi)}{2},\quad p_{1,-1}(\xi)=\frac{1- \beta_1(\xi)}{2},  \quad \xi \in \R^{d-1}.$$
Let $\{x^n\}_{n \in \N}$ be a sequence in $\ZZ^d$ such that $\bar x^n = x^n/\sqrt{n}$ converges to $x$.
For each $n$ define a Markov chain $\{X^n_k\}_{k \in \N_0}$ such that $X^n_0 = x^n$ and for $k \ge 0$, $v \in \ZZ^d$
\begin{align}
	P(\Del X_{k}^n &= u \mid X_k^n = v) = \frac{1}{2^d}, \; u \in \{-1,1\}^d, \mbox{ if } v_1 \neq 0\nonumber\\
	P(\Del X_{k}^n &= \sum_{i=2}^d e_i \bar \beta_i(v'/\sqrt{n}) + u \mid X_k^n = v)
\nonumber
\\
&\qquad=\prod_{i=1}^d p_{i, u_i}(v'/\sqrt{n}), \; u = (u_1, \ldots, u_d) \in \{-1,1\}^d, \mbox{ if } v_1 = 0.\nonumber \\
	\label{eq:delx}
	\end{align}
	Here $v = (v_1, v') \in \ZZ \times \ZZ^{d-1}$, $\{e_i\}_{i=1}^d$ is the standard coordinate basis in $\R^d$ and $\Del Z_k = Z_{k+1} - Z_k$ for an $\R^d$-valued sequence $\{Z_k\}$.
	Since the transition probabilities of this chain agree with those of a homogeneous
    random walk (RW)
	away from $\{v_1=0\}$ we can couple this chain with a $d$-dimensional RW $\{W^n_k\}_{k \in \N_0}$ such that the following identity holds:
	\begin{equation}
	X^n_k = x^n + \sum_{i=0}^{k-1} \left(\Del W^n_i 1_{U^n_i \neq 0} + \Del X^n_i 1_{U^n_i = 0}\right)
	\label{eq:stateeq}
\end{equation}
	where we write $X^n_i = (U^n_i, Y^n_i) \in \R \times \R^{d-1}$.  Letting $\clf^{n}_k = \sigma\{X^n_j, W^n_j, j \le k\}$ we have
	$$
	E(\Del X^n_{k} \mid \clf^n_k ) = \beta(Y^n_k/\sqrt{n}) 1_{U^n_k = 0}.$$
	We can therefore write
	\begin{equation}
		\label{eq:08}
		X^n_k = x^n + W^n_k + \sum_{i=0}^{k-1} \left(\beta(Y^n_i/\sqrt{n}) + M^n_i\right)1_{U^n_i = 0}\end{equation}
	where
	$|M^n_i| \le (3+ 2 \sup_{\xi}|\beta(\xi)|)\sqrt{d} \equiv  c_1$
    and $\hat M^n_k = \sum_{i=0}^{k-1} M^n_i 1_{U^n_i = 0}$ is a martingale.
	
	Define $L^n_0 =0$ and for $k \ge 1$, $L^n_k = \sum_{i=0}^{k-1} 1_{U^n_i = 0}$.
	Note that the increasing process associated with $\hat M^n$ is bounded by
	$2c_1^2L^n$,  and thus by Burkholder's inequality,
	\begin{equation}\label{09}
	E[\max_{i\le k}(\hat M^n_i)^2]\le c_2 E[L^n_k]
	\end{equation}
where $c_2 = 8c_1^2$.
Next note that
$$\Del |U^n_k| = \Del U_k^n 1_{U_k^n > 0} - \Del U_k^n 1_{U_k^n < 0} + 1_{U_k^n = 0}.$$
Thus
\begin{equation}\label{19}|U_k^n| = |x_1^n| + Z_k^n + L_k^n,\end{equation}
where $Z_k^n = \sum_{i=0}^{k-1} \sgn(U^n_i)\Del U^n_i$.
Recalling \eqref{eq:delx} and that $\sgn(0) =0$ we see that
\begin{equation}
\label{eq:ezeq0}
E(Z_k^n) = E\left(\sum_{i=0}^{k-1} 1_{U_i^n \neq 0}\, \sgn(U^n_i) E(\Del U^n_i\mid \clf^n_i)\right) = 0.
\end{equation}
Next let $\{\zeta_i\}$ be an i.i.d. sequence of $\{-1,1\}$-valued random variables, independent of all the sequences introduced above, such that $P(\zeta_i = 1) = P(\zeta_i = -1) = 1/2$
and let
$$Z_k^{*,n} = Z_k^n + \sum_{i=0}^{k-1} \zeta_i 1_{U_i^n = 0}.$$
By construction $\{Z^{*,n}_k\}$ is a simple random walk on $\ZZ$.  Also, as for \eqref{09}
\begin{equation}
  \label{11}
  E[\max_{i\le k}(Z^{*,n}_i-Z^n_i)^2]\le 8E[L^n_k].
\end{equation}
 Define now a continuous time process
 \[
 \bar W^n(t)=\frac{1}{\sqrt n}W^n_{\lf nt\rf},\qquad t\in[0,\iy)
 \]
and use the same rescaling to define continuous time versions
 $\bar X^n,\bar U^n,\bar Y^n,\bar L^n,\bar Z^n,\bar Z^{*,n}$ of the discrete time processes
 $X^n,U^n,Y^n,L^n,Z^n,Z^{*,n}$.
Then, from \eqref{eq:08},
\begin{equation}
	\label{13}
	\bar X^n(t) = \bar x^n + \bar W^n(t) + \int_{[0,\frac{\lf nt\rf}{n}]} \beta (\bar Y^n(s-)) d\bar L^n(s) + \eps^n(t),
\end{equation}
where by \eqref{09}
\begin{equation}
	\label{eq:martbd}
E[\max_{s\le t}(\eps^{n}(s))^2]\le \frac{c_2}{\sqrt n}E[\bar L^n(t)],
\end{equation}
By \eqref{19} and \eqref{eq:ezeq0} we have
$$
E(\bar L^n(t)) \le E|U^n(t)|
\le |\bar x_1^n| + \sqrt{t},$$
where the last inequality follows on observing that $|U^n|$ has the same law as the absolute value of a RW starting from $\bar x^n_1$.
Thus
\begin{equation}\label{eq:25}
\sup_n E(\bar L^n(t)) < \infty \mbox{ for all } t \ge 0.
\end{equation}
Using this in \eqref{eq:martbd} we get that $\eps^n$ converges to $0$ in probability, uniformly on compacts (u.o.c.).
Also, recalling that $\{Z^{*,n}_k\}$ is a simple random walk, we have by Donsker's theorem that
$\bar Z^{*,n}$ converges in distribution in $\DD_{\R}$ to a standard Brownian motion $Z$.  From \eqref{11} and \eqref{eq:25}
we have that $\bar Z^{*,n} - \bar Z^n$ converges to $0$ in probability, u.o.c. and so $\bar Z^n$ converges to $Z$ weakly in $\DD_{\R}$ as well.

Consider the 1-dimensional Skorohod map $\Gam:  \DD_{\R} \to \DD_{\R_+}$,
$$\Gam(g)(t) = g(t) - \inf_{0\le s \le t}(g(s) \wedge 0), t \ge 0, \; g \in \DD_{\R}$$
and recall that $-\inf_{0\le s \le t}(g(s) \wedge 0)$ is the unique nonnegative, nondecreasing RCLL function $h$ that satisfies  $f(t) = g(t) + h(t) \ge 0$ for all $t \ge 0$ and
$\int_{[0,\infty)} 1_{(0,\infty)}(f(t)) dh(t) = 0$.
Writing
$$|\bar U^n(t)| = \left(|\bar x^n_1| + \bar Z^n(t) - \frac{1}{\sqrt n} 1_{\bar U^n(t)=0}\right) + \bar L^n\left( \frac{\lf nt\rf + 1}{n}\right),$$
we see that
$$|\bar U^n(t)| = \Gam (|\bar x^n_1| + \tilde Z^n)(t),\quad \mbox{where}\quad \tilde Z^n(t) = \bar Z^n(t) - \frac{1}{\sqrt n} 1_{\bar U^n(t)=0}.$$
Since $|\tilde Z^n(t) - \bar Z^n(t)| \le 1/\sqrt{n}$, we have that $\tilde Z^n$ converges to $Z$ weakly in $\DD_{\R}$.  By the continuity of $\Gam$, we now have that
$(\bar Z^n,|\bar U^n|,\bar L^n)$ converges in distribution in $\DD_{\R^3}$ to
$(Z,R,L)$,
where
\[
R=|x_1|+Z-\min_{s\le\cdot}[(|x_1|+Z(s))\wedge 0], \qquad L= -\min_{s\le\cdot}[(|x_1|+Z(s))\wedge 0].
\]

 Next, using again Donsker's theorem, $\bar W^n$ converges in distribution to a $d$-dimensional standard Brownian motion. In particular,
 $\bar W^n$ is $C$-tight.  Combining this with \eqref{13} and the weak convergence of $\eps^n$ to $0$ and $\bar L^n$ to $L$ we see that
$\bar X^n$ is $C$-tight as well.
 Denote by $(W,X,U,Y,Z,R,L)$ a subsequential weak limit of $(\bar W^n, \bar X^n,\bar U^n, \bar Y^n, \bar Z^n, |\bar U^n|, \bar L^n)$.
The integral  in \eqref{13} converges to $\int_0^\cdot \beta(Y(s))dL(s)$  by Theorem 2.2 of \cite{kur-pro}.
We therefore obtain
\[
X(t)=x+W(t) +\int_0^t \beta(Y(s))dL(s),\qquad \qquad t\ge0.
\]
Since $R$ is the limit of $|\bar U^n|$, the relation $R=|U|$ must hold.  Also $X = (U, Y)$.
For a continuous and bounded function $G: \bar C \times \CC_{\R} \to \R$ and $0 \le t \le s < \infty$
$$
\left | E \left (G\left(\bar W^n(\cdot \wedge t), \bar X^n(\cdot \wedge t), \bar L^n(\cdot \wedge t), \bar Z^n(\cdot \wedge t) \right) \left[ \bar W^n(s) - \bar W^n(t)\right] \right)\right|
\le \frac{2 \|G\|_{\infty}}{\sqrt n}$$
and thus
$$
 E \left (G\left(\bar W(\cdot \wedge t), \bar X(\cdot \wedge t), \bar L(\cdot \wedge t), \bar Z(\cdot \wedge t) \right) \left[ \bar W(s) - \bar W(t)\right] \right)
=0.$$
Consequently $W$ is a standard $d$-dimensional $\clf_t$-Brownian motion where
$\{\clf_t = \sigma\{X(\cdot \wedge t),W(\cdot \wedge t), L(\cdot \wedge t), Z(\cdot \wedge t)\}\}_{t \ge 0}$.
Similarly $Z$ is a 1-dimensional $\clf_t$-Brownian motion.
Next note that
\begin{equation*}
	|U(t)| = |x_1| + Z(t) + L(t)
\end{equation*}
and by Tanaka's formula, since $L$ increases only when $|U| = 0$,
$$
|U(t)| = |x_1| + \int_0^t \sgn(U(s)) dW_1(s) + L^U(t), $$
where $L^U$ is the symmetric local time of $U$ at $0$.  Combining the above two displays we have
$L = L^U$.  Thus the processes $(X,W,L)$ satisfy all the properties listed at the beginning of the subsection and consequently we have proved the existence of a weak solution.
\qed
\section{Brownian particles with skew-elastic collisions}\label{sec4}
\beginsec

As an application we consider a model introduced by Fernholz et al.\ \cite{FIK}
for the dynamics of
a pair of 1-dimensional Brownian particles $X_1$ and $X_2$ that exhibit various possible
types of interaction when they collide. The equations involve the local time at zero of the relative
position, and the types of interaction are determined
by the coefficients in front of the local time terms.
For a continuous real semimartingale $Z$, let $L^Z_+$ be defined by the relation
$$L^Z_+(t) = Z^+(t) - Z^+(0) - \int_0^t 1_{Z(s) > 0} dZ(s)$$
and let $L^Z_{-} = L^{-Z}_+$.  Then the symmetric local time of $Z$ at $0$, $L^Z$,
is given as $L^Z = L^Z_+ + L^Z_{-}$.
The motion of the particles is described by the set of equations
\begin{align}
	\label{eq:skewcol1}
	dX_1(t) &= k_1(X(t)) dt + \frac{1}{\sqrt{2}} dB_1(t)
	+ \frac{1 - \zeta_1(X(t))}{2} dL_+^{X_1-X_2}(t)
	+ \frac{1 - \eta_1(X(t))}{2} dL_+^{X_2-X_1}(t)\nonumber\\
	dX_2(t) &= k_2(X(t)) dt + \frac{1}{\sqrt{2}} dB_2(t)
	+ \frac{1 - \zeta_2(X(t))}{2} dL_+^{X_1-X_2}(t)
	+ \frac{1 - \eta_2(X(t))}{2} dL_+^{X_2-X_1}(t).\nonumber \\
\end{align}
Here  $(B_1,B_2)$ is a planar Brownian motion,
$X = (X_1, X_2)$ and  $k_i, \zeta_i, \eta_i$, $i=1,2$, are bounded measurable functions from $\R^2$ to $\R$
(additional conditions on  $\zeta_i, \eta_i$, $i=1,2$ will be introduced below).
In \cite{FIK} the coefficients $\zeta_i$ and $\eta_i$ are assumed to be constant rather than state-dependent.
Different types of behavior are obtained by different choices of these constants.
For example, when $\zeta_1=\eta_1=1$ (resp., $\zeta_2=\eta_2=1$),
the local times disappear from the equation for $X_1$ (resp., $X_2$) and
so the collisions do not affect the first (resp., second) particle.
When both conditions hold, the motions are completely frictionless.
On the other extreme, when $\zeta_2=\eta_2=1$ and $\zeta_1=\eta_1=-1$,
the trajectory of $X_1$ bounces off
that of $X_2$ as if it were a perfectly reflecting boundary.
See \cite{FIK} for more general conditions under which frictionless motion and perfect
reflection are attained.
Other combinations of the constants give rise to a whole range of elastic collisions.

The equations presented here have state-dependent local time coefficients.
This allows one to
model variability in the type of collision, where the type is determined by the collision position.
(The paper \cite{FIK} also emphasizes rank-dependent motion, and to this end the
drift and diffusion coefficients considered are constant on each of the sets $\{X_1\le X_2\}$ and $\{X_1>X_2\}$.
Such a drift can be realized in our model by selecting $k_i$ suitably.
However, we only consider constant diffusion coefficient here.)
For the equations with constant local time coefficients,
\cite{FIK} obtain strong existence and pathwise uniqueness. Our goal here is to exhibit
weak existence and uniqueness for equations \eqref{eq:skewcol1}.

For $x\in\R^2$ let
$$\zeta(x) = 1 + \frac{\zeta_1(x) - \zeta_2(x)}{2}, \quad \eta(x) = 1 - \frac{\eta_1(x) - \eta_2(x)}{2},
$$
$$
\bar \zeta(x) = 1 - \frac{\zeta_1(x)+\zeta_2(x)}{2},\quad
\bar \eta(x) = 1 - \frac{\eta_1(x)+\eta_2(x)}{2}.
$$
We assume that for all $x \in \R^2$,
$$\zeta(x) \ge 0, \quad \eta(x) \ge 0, \quad \eta(x) + \zeta(x) \neq 0, \quad \alpha(x) = \frac{\eta(x)}{\eta(x) + \zeta(x)} \in [0,1].$$
Consider the equations associated with the linear transformations
$$Y = X_1 - X_2, \quad U = X_1 + X_2.$$
If $X = (X_1, X_2)$ solves  \eqref{eq:skewcol1} then $(Y,U)$ will solve
\begin{align}
	\label{eq:skewcol2}
	dY(t) &= (k_1(X(t)) -k_2(X(t))) dt + dW_1(t)
+ (1 - \zeta(X(t))) dL_+^Y(t) + (\eta(X(t))-1) dL_-^Y(t)\nonumber\\
	dU(t) &= (k_1(X(t)) + k_2(X(t))) dt + dW_2(t) +  \bar \zeta(X(t)) dL_+^Y(t) +  \bar \eta(X(t)) dL_-^Y(t),\nonumber\\
\end{align}
where $W_1 = \frac{1}{\sqrt{2}}(B_1 - B_2)$, $W_2 = \frac{1}{\sqrt{2}}(B_1 + B_2)$.
Using the relations
$$L^Y_+(t) - L^Y_-(t) = \int_0^t 1_{[Y(s)=0]} dY(s), \quad L^Y(t) =  (L^Y_+(t) + L^Y_{-}(t)), \; t \ge 0,$$
and $\int_{[0, \infty)} 1_{[Y(s)=0]} ds = 0$, we see that
$$
L_+^Y(t) =  \int_0^t \alpha(X(s)) dL^Y(s), \quad L_-^Y(t) =  \int_0^t (1-\alpha(X(s))) dL^Y(s).$$
Define $\psi: \R^2 \to \R^2$ as $\psi(u,y) = (\frac{u+y}{2}, \frac{u-y}{2})$.  Then \eqref{eq:skewcol2} can be rewritten as
\begin{align}
	\label{eq:skewcol3}
	dY(t) &= a_1(U(t),Y(t)) dt + dW_1(t) + b_1(U(t)) dL^Y(t)\nonumber\\
	dU(t) &= a_2(U(t),Y(t)) dt + dW_2(t) + b_2(U(t)) dL^Y(t),\nonumber\\
\end{align}
where
\begin{align*}
	a_1(u,y) &= (k_1-k_2)\circ\psi(u,y), \quad a_2(u,y) = (k_1+k_2)\circ \psi(u,y),\\
b_1(u) & = \beta_1 \circ \psi(u,0), \quad b_2(u) = \beta_2 \circ \psi(u,0),\\
\beta_1(x) &= (2\alpha(x)-1),
\quad \beta_2(x) = \bar \zeta(x)\alpha(x) +  \bar \eta(x)(1-\alpha(x)).
\end{align*}
The following is now an immediate consequence of Theorem \ref{thm:mainthm}.
\begin{theorem}
	\label{thm:skewcoll}
	Suppose that $\zeta_i, \eta_i$, $i = 1,2$ and $\alpha$ are Lipschitz. Then equations \eqref{eq:skewcol3} and \eqref{eq:skewcol1}, with any given initial conditions, have a unique weak solution.
\end{theorem}

\noi{\bf Proof.}
Note our assumptions ensure that $b$ is Lipschitz with $b_1$ taking values in $[-1,1]$.
Thus existence of a unique weak solution of \eqref{eq:skewcol3} is an immediate consequence of
Theorem \ref{thm:mainthm}.
Unique solvability of  \eqref{eq:skewcol1} follows on observing that $(X_1, X_2)$ solves \eqref{eq:skewcol1} with driving Brownian motions $(B_1, B_2)$ if and only if
$(Y,U)$ solves \eqref{eq:skewcol3} with driving Brownian motions $(W_1,W_2)$, where $Y = X_1-X_2$, $U = Y_1 + Y_2$, $W_1 = \frac{1}{\sqrt{2}}(B_1 - B_2)$, $W_2 = \frac{1}{\sqrt{2}}(B_1 + B_2)$.
\qed


%

{\sc

\bigskip\noi
R. Atar\\
Department of Electrical Engineering\\
Technion\\
Haifa 32000, Israel\\
email: atar@ee.technion.ac.il

\skp

\noi
A. Budhiraja\\
Department of Statistics and Operations Research\\
University of North Carolina\\
Chapel Hill, NC 27599, USA\\
email: budhiraj@email.unc.edu

}

\end{document}